# Investigation on Bifurcation of Map Satisfying Condition of Periodic Doubling Bifurcation

Dr. Vice-prof. *Kim Sang Mun*

College of Mathematics, **Kim Il Sung** University

**Abstract** In this paper, we consider sufficient conditions for an invariant double circle to occur in a 1 parameter discrete dynamical systems on R×S.

We consider the sufficient conditions for invariant double circle to occur in a discrete dynamical system.

In the previous papers, were investigated sufficient conditions for periodic doubling bifurcation to occur, in 1 parameter-1 dimension systems and sufficient conditions for Naimark-Sacar bifurcation to occur, in 1 parameter-2 dimension systems.

Then, was investigated a sufficient condition for periodic doubling bifurcation to occur in the case that higher order derivative is equal to zero, in 1 parameter-1 dimension systems.

And were studied sufficient conditions for Hopf- bifurcation to occur in the case that partial derivative of higher order is equal to zero and the normal form of Hopf- bifurcation. In this paper we consider conditios for invariant double circle to occur for the discrete dynamical system on R×S.

**Theorem 1** If $\lambda$ $(3 < \lambda < \lambda_0)$ irrational number and $F_\lambda(r, [\theta]) = (f_\lambda(r), H([\theta]))$ then there exists $\lambda_0$ $(>3)$ such that the

Let $f_\lambda(r) = \lambda r(1-r)$, $r > 0$, $\lambda > 0$, $H([\theta]) = [\theta + \alpha]$, $[\theta] \in S^1$, here $\alpha$ is an irrational number; $F_\lambda(r, [\theta]) = (f_\lambda(r), H([\theta]))$ (where [ ] is a symbol of equivalent class of $\mathbf{R}/\mathbf{Z}$).

Then there exists $\lambda_0$ $(>3)$ such that for each

$\lambda$ $(3 < \lambda < \lambda_0)$ there exist the invariant circles $\Gamma_\lambda^1$, $\Gamma_\lambda^2$ with respect to $F_\lambda^2$ and the following are true:

① $\Gamma_\lambda^1 \cap \Gamma_\lambda^2 = \phi$

② $F_\lambda(\Gamma_\lambda^1) = \Gamma_\lambda^2$, $F_\lambda(\Gamma_\lambda^2) = \Gamma_\lambda^1$

③ $F_\lambda(\Gamma_\lambda^1 \cup \Gamma_\lambda^2) = \Gamma_\lambda^1 \cup \Gamma_\lambda^2$

**Proof** Since $f_\lambda(r)$ satisfies the condition of periodic doubling bifurcation, there exists some $\lambda_0$ $(>3)$ such that for each $\lambda$ $(3 < \lambda < \lambda_0)$, there exist $r_1(\lambda), r_2(\lambda) > 0$ satisfying $f_\lambda(r_1) = r_2$, $f_\lambda(r_2) = r_1$, $r_1 \neq r_2$.[1]

Let define $\Gamma_\lambda^1$ and $\Gamma_\lambda^2$ as follows:

$$\Gamma_\lambda^1 := \{r_1(\lambda)\} \times S^1, \quad \Gamma_\lambda^2 := \{r_2(\lambda)\} \times S^1$$

Let write $r_1(\lambda)$, $r_2(\lambda)$ as $r_1, r_2$ and suppose $(r_1, [0]) \in \Gamma_\lambda^1$, then

$$\left.\begin{array}{l} F_\lambda^1(r_1, [0]) = (r_2, [\alpha]) \in \Gamma_\lambda^2 \\ F_\lambda^2(r_1, [0]) = (r_1, [2\alpha]) \in \Gamma_\lambda^1 \\ F_\lambda^3(r_1, [0]) = (r_2, [3\alpha]) \in \Gamma_\lambda^2 \\ F_\lambda^4(r_1, [0]) = (r_1, [4\alpha]) \in \Gamma_\lambda^1 \\ \vdots \end{array}\right\}$$

If $k$ is a nonnegative integer, then

$$F_\lambda^{2k}(r_1, [0]) = (r_1, [2k\alpha]) \in \Gamma_\lambda^1, \quad F_\lambda^{2k+1}(r_1, [0]) = (r_2, [(2k+1)\alpha]) \in \Gamma_\lambda^2.$$

First prove that $\{F_\lambda^{2k}(r_1, [0])\}$ is dense in $\Gamma_\lambda^1$.

The set $\{F_\lambda^{2k}(r_1, [0])\}$ consists of different points. Therefore, since $\Gamma_\lambda^1$ is compact, this sequence has an accumulation point in $\Gamma_\lambda^1$. And $\{r_1\} \times S^1$ is a complete metric space with distance.

$$d((r_1, [x]), (r_1, [y])) = \inf\{|x - y + k|\,;\, k \in \mathbf{Z}\}$$

Since $\{F_\lambda^{2k}(r_1, [0])\}$ has an accumulation point in $\Gamma_\lambda^1$

$$\forall \varepsilon > 0,\ \exists k_1, k_2\ (k_1 > k_2) \in \mathbf{Z}^+;\ d(F_\lambda^{2k_1}(r_1, [0]), F_\lambda^{2k_2}(r_1, [0])) = \inf\{|(2k_1 - 2k_2)\alpha + k|\,;\, k \in \mathbf{Z}\} < \varepsilon.$$

If $p := |2k_1 - 2k_2|$ then p is an even number.

Therefore $d(F_\lambda^p(r_1, [0]), (r_1, [0])) = d(F_\lambda^{2k_1}(r_1, [0]), F_\lambda^{2k_2}(r_1, [0])) < \varepsilon$.

If $k' \in N$ then $d(F_\lambda^{k'p}(r_1, [0]), F_\lambda^{(k'-1)p}(r_1, [0])) = d(F_\lambda^p(r_1, [0]), (r_1, [0])) < \varepsilon$.

Therefore $\{r_1\} \times S^1$ is divided into parts whose length is smaller than $\varepsilon$ by the sequence

$$F_\lambda^p(r_1, [0]),\ F_\lambda^{2p}(r_1, [0]),\ \cdots.$$

Therefore $\{F_\lambda^{2k}(r_1, [0])\}$ is dense in $\Gamma_\lambda^1$. That is

$$\overline{\{F_\lambda^{2k}(r_1, [0])\}} = \Gamma_\lambda^1. \tag{1}$$

Similarly,

$$\overline{\{F_\lambda^{2k+1}(r_1, [0])\}} = \Gamma_\lambda^2. \tag{2}$$

From (1), (2) and the property of invariant sets, $\Gamma_\lambda^1 = F_\lambda^2(\Gamma_\lambda^1)$, $\Gamma_\lambda^2 = F_\lambda^2(\Gamma_\lambda^2)$ and from the definition $\Gamma_\lambda^1 \cap \Gamma_\lambda^2 = \phi$. Thus the result ① is proved.

Let prove ②. If $p \in F_\lambda(\Gamma_\lambda^1)$ then there exists some $u \in \Gamma_\lambda^1$ such that $p = F_\lambda(u)$.

And there exists $\{(r_1, [2k_1\alpha])\} \subset \Gamma_\lambda^1$ such that $u = \lim_{k_1 \to \infty}(r_1, [2k_1\alpha])$.

$\overline{\varphi}: \mathbf{R} \to \mathbf{R}/\mathbf{Z}$, $\overline{\varphi}(x) := [x]$, $\varphi: \mathbf{R} \to S^1$, $\varphi(x) := e^{2\pi i x}$, $\pi: \mathbf{R}/\mathbf{Z} \to S^1$, $\pi[x] = e^{2\pi i x}$, $x \in \mathbf{R}$

then $\pi$ is a homeomorphism and

$$\mathbf{R} \xrightarrow{\overline{\varphi}} \mathbf{R}/Z \xrightarrow{\pi} S \xleftarrow{\varphi} \mathbf{R}$$

$\overline{\varphi}(x) = \pi \circ \varphi(x)$, therefore $\overline{\varphi}$ is continuous.

Since $F_\lambda(r, [\theta]) = (f_\lambda(r), H([\theta])) = (\lambda r(1-r), \overline{\varphi}(\theta + \alpha))$, $F_\lambda$ is continuous with respect to $(r, [\theta])$. Thus, $p = F_\lambda(u) = \lim_{k_n \to \infty}(r_2, [(2k_1+1)\alpha])$ and since $\Gamma_\lambda^2$ is compact, $p = F_\lambda(u) \in \Gamma_\lambda^2$.

Therefore
$$F_\lambda(\Gamma_\lambda^1) \subset \Gamma_\lambda^2 \tag{3}$$

If $v \in \Gamma_\lambda^2$ then there exist the sequence $\{r_2, [(2k_n+1)\alpha]\} \subset \Gamma_\lambda^2$ such that
$$v = \lim_{k_n \to \infty}(r_2, [(2k_1+1)\alpha]).$$

Since $F_\lambda$ is continuous, $v = F_\lambda(\lim_{k_n \to \infty}(r_1, [2k_n\alpha]))$ and by the compactness of $\Gamma_\lambda^1$.

If $u = \lim_{k_n \to \infty}(r_1, [2k_n\alpha])$ then $v = F_\lambda(u)$, $u \in \Gamma_\lambda^1$. Thus
$$\Gamma_\lambda^2 \subset F_\lambda(\Gamma_\lambda^1). \tag{4}$$

By (3), (4), $F_\lambda(\Gamma_\lambda^1) = \Gamma_\lambda^2$ and $F_\lambda(\Gamma_\lambda^2) = \Gamma_\lambda^1$.

The proof of ③ is omitted.

**Theorem 2** Assume that $f: \mathbf{R}^2 \to \mathbf{R}$ satisfies the condition of periodic doubling bifurcation in a neighborhood of $f(\lambda, \cdot)$ and $g: \mathbf{R}^2 \ni (\lambda, r) \mapsto g(\lambda, r) \in \mathbf{R}$ is continuous with respect to $r$ and for 2-periodic points $r_1$, $r_2$ of $f(\lambda, \cdot)$, $g(\lambda, r_1)$, $g(\lambda, r_2)$ are irrational numbers, and $\Theta_\lambda: \mathbf{R} \times S^1 \to S^1$, $\Theta_\lambda(r, [\theta]) := [\theta] + [g(\lambda, r)]$.

If $F_\lambda(r, [\theta]) := (f(\lambda, r), \Theta_\lambda(r, [\theta]))$ then there exists $\lambda_1 \in \mathbf{R}$ such that for each $\lambda \in (\lambda_0, \lambda_1)$ (or $\lambda \in (\lambda_1, \lambda_0)$) there exist invariant circle $\Gamma_\lambda^1$, $\Gamma_\lambda^2$ with respect to $F_\lambda^2$ such that

① $\Gamma_\lambda^1 \cap \Gamma_\lambda^2 = \phi$, $F_\lambda^2(\Gamma_\lambda^1) = \Gamma_\lambda^2$, $F_\lambda^2(\Gamma_\lambda^1) = \Gamma_\lambda^2$

② $F_\lambda(\Gamma_\lambda^1) = \Gamma_\lambda^2$, $F_\lambda(\Gamma_\lambda^2) = \Gamma_\lambda^1$

③ $F_\lambda(\Gamma_\lambda^1 \cup \Gamma_\lambda^2) = \Gamma_\lambda^1 \cup \Gamma_\lambda^2$